\documentclass[12pt]{article}
\usepackage{amssymb,amsmath,amsthm}

\newtheorem{thm}{Theorem}
\newtheorem*{thma}{Theorem A}
\newtheorem*{thmb}{Theorem B}
\newtheorem*{thmc}{Theorem C}

\newtheorem{defn}{Definition}
\newtheorem{que}{Question}
\newtheorem{rem}{Remark}

\newtheorem{cor}{Corollary}

\newcommand {\N}  {{\mathbb N}}

\newcommand {\R}  {{\mathbb R}}

\begin{document}
\title{Higher Order Oscillation and Uniform Distribution~~}

\date{}
\author{Shigeki AKIYAMA\footnote{The first author is supported by Japanese Society for the Promotion of Science (JSPS), Grant in aid 26287017.}
 and Yunping JIANG~\footnote{The second author is partially supported by the collaboration grant from the Simons Foundation [grant number 199837] and the CUNY collaborative incentive research grants [grant number 2013] and awards from PSC-CUNY and grants from NSFC [grant numbers 11171121 and 11571122].}}
\maketitle

\abstract{It is known that the M\"obius function in number theory is higher order oscillating. 
In this paper we show that there is another kind of higher order oscillating sequences 
in the form $(e^{2\pi i \alpha \beta^{n}g(\beta)})_{n\in \N}$, for a non-decreasing
twice differentiable function $g$ with a mild condition.
This follows the result we prove in this paper that
for a fixed non-zero real number $\alpha$ and almost all real numbers $\beta >1$ (alternatively, for a fixed real number $\beta >1$ and almost all real numbers $\alpha$) and for all real polynomials $Q(x)$, sequences $\big(\alpha \beta^{n}g(\beta)+Q(n)\big)_{n\in \N}$ are uniformly distributed modulo $1$. 
}

\section{Introduction}
We denote by $\N$ the set of positive integers. Suppose ${\bf c}=(c_{n})_{n\in \N}$, is a sequence of complex numbers. In~\cite{FanJiang}, an oscillating sequence is defined for the purpose of the study of Sarnak's conjecture (see~\cite{Sa1,Sa2}) which is stated as 
the M\"obius function is linearly disjoint from all zero entropy flows. Let us recall the definition of an oscillating sequence.

\medskip
\begin{defn}[Oscillation]~\label{os}
The sequence ${\bf c}=(c_{n})_{n\in \N}$ is said to be oscillating if 
\begin{equation}~\label{control}
\sum_{n=1}^{N} |c_{n}|^{\lambda} =O(N) \hbox{ for some $\lambda>1$}.
\end{equation}  
and if
\begin{equation}~\label{eqhoosw}
\lim_{N\to \infty} \frac{1}{N} \sum_{n=1}^{N} c_{n}e^{2\pi i n t}=0, \;\; \forall \; 0\leq t<1.
\end{equation} 
\end{defn}

Recall that the M\"obius function $\mu (n)$ is, by definition, $\mu(n)=1$ if $n=1$; $\mu(n)=(-1)^{r}$ if $n=p_{1}\cdots p_{r}$ for $r$ distinct prime numbers $p_{i}$; $\mu (n)=0$ if $p^{2}|n$ for some prime number $p$. 
The M\"obius sequence ${\bf u}=(\mu (n))_{n\in \N}$ is the one generated by the M\"obius function. Due to Davenport's theorem~\cite{Da}, the M\"obius sequence is oscillating. 

We proved in~\cite{FanJiang} that any oscillating sequence is linearly disjoint from all minimally mean attractable (MMA) and minimally 
mean-L-stable (MMLS) flows.   
In the same paper, 
we further proved that flows defined by all $p$-adic polynomials of integral coefficients,  all $p$-adic rational maps with good reduction, 
all automorphisms of the $2$-torus with zero topological entropy, all diagonalizable affine maps of the $2$-torus with zero topological entropy, 
all Feigenbaum maps, and all  orientation-preserving circle homeomorphisms are MMA and MMLS. 
Therefore, we confirmed Sarnak's conjecture for these flows which form a large class of zero topological entropy flows. 
However, it is also shown in~\cite[Example 7]{FanJiang}, 
only the oscillation property is not enough for the study of Sarnak's conjecture. We need a higher order oscillation condition in the study of Sarnak's conjecture. We have two versions of a definition of a higher order of oscillating sequence as appeared in~\cite{J} (see also~\cite[Remark 8]{FanJiang}). 

\medskip
\begin{defn}[Weaker Version of Higher Order Oscillation]~\label{hoosw}
We call the sequence 
${\bf c}=(c_{n})_{n\in \N}$ 
a higher order oscillating sequence of order $m\geq 2$ if it satisfies (\ref{control}) and if 
\begin{equation}~\label{hoosweq}
\lim_{N\to \infty}\frac{1}{N} \sum_{n=1}^{N} c_{n}e^{2\pi i n^{k}t}=0, \;\;\forall\; 1\leq k\leq m,\;\;\forall\; 0\leq t\leq 1.
\end{equation}
\end{defn}

Thanks to Hua's result~\cite{Hua}, we knew that the M\"obius sequence ${\bf u}$ is a higher order oscillating sequence of order $m$ for all $m\geq 2$ in this weaker version of the definition (Definition~\ref{hoosw}). Actually, according to~\cite{LZ}, we have that for any $A>0$,
\begin{equation}~\label{hoosd}
\sum_{n=1}^{N} c_{n}e^{2\pi i n^{k}t}=O_{A}\Big( N (\log N)^{-A}\Big), \;\;\forall\; 1\leq k\leq m,\;\;\forall\; 0\leq t\leq 1.
\end{equation}

\medskip
\begin{defn}[Stronger Version of Higher Order Oscillation]~\label{hooss}
We call the sequence 
$
{\bf c}=(c_{n})_{n\in \N}
$ 
a higher order oscillating sequence of order $m\geq 2$ if it satisfies (\ref{control}) and if 
\begin{equation}~\label{hoosseq}
\lim_{N\to \infty} \frac{1}{N} \sum_{n=1}^{N} c_{n}e^{2\pi i P(n)}=0
\end{equation} 
for every real polynomial $P$ of degree $\leq m$. 
\end{defn}

We prove in~\cite{J} that any higher order oscillating sequence of order $d$ is linearly disjoint 
from all affine distal flows on the $d$-torus for all $d\geq 2$. One consequence of this result is that 
any higher order oscillating sequence of order $2$ is linearly disjoint 
from all affine flows on the $2$-torus with zero topological entropy. 
In their paper~\cite[Lemma 2.1]{LS}, Liu and Sarnak showed that the M\"obius sequence ${\bf u}$
is also an higher order oscillating sequence of order $m$ for all $m\geq 2$ in the stronger version of the definitin (Definition~\ref{hooss}). They actually showed an estimation like the one in (\ref{hoosd}). 
Combining this with our main result in~\cite{J}, we reconfirms Sarnak's conjecture for all affine flows 
on the $2$-torus with zero topological entropy and 
for all affine distal flows on the $d$-torus for all $d>2$. 
Then we have the following interesting question.

\medskip
\begin{que}~\label{que}
Is there another kind of higher order oscillation sequences either in the weaker version of the definition or in the stronger version of the definition except for the one generated by an arithmetic function? 
\end{que}

We study this question in this paper. 

\section{Statement of the Main Result}

For a real number $x$, let $[x]$ denote the integer part of $x$, that is, the greatest integer $\leq x$; let 
$$
\{x\}=x-[x]
$$ 
be the fractional part of $x$, or the residue of $x$ modulo $1$. 
 
\medskip
\begin{defn}[Uniform Distribution]
We say a sequence ${\bf x}=(x_{n})_{n\in \N}$ of real numbers is uniformly distributed modulo $1$ (abbreviated u.\,d.\,mod\;$1$)  if for any $0\leq a<b\leq 1$,
we have
$$
\lim_{N\to \infty} \frac{\#(\{ n\in [1,N] \;|\; \{x_{n}\} \in [a,b]\})}{N}=b-a.
$$
 \end{defn}

We state three results in the u.\,d.\,mod\;$1$ theory. 

\medskip
\begin{thma}[The Weyl Criterion]~\label{wc}
The sequence ${\bf x}=(x_{n})_{n\in \N}$ is u.\,d.\,mod\;$1$ if and only if 
$$
\lim_{N\to \infty} \frac{1}{N} \sum_{n=1}^{N} e^{2\pi i h x_{n}} =0 \quad \hbox{for all integers $h\not=0$}.
$$
\end{thma}

\medskip
\begin{thmb}[Koksma's Theorem]~\label{kth}
Let $(y_{n}(x))_{n\in \N}$ be a sequence of real valued $C^{1}$ functions defined on an interval $[a,b]$. 
Suppose $y_{m}'(x)-y_{n}'(x)$ is monotone on $[a,b]$ for any two integers $m\not=n$ 
and suppose 
$$
\inf_{m\not=n}\min_{x\in [a,b]} |y_{m}'(x)-y_{n}'(x)| >0.
$$
Then for almost all $x\in [a,b]$, the sequence ${\bf y} = (y_{n}(x))_{n\in \N}$ is u.\,d.\,mod\;$1$.
\end{thmb}

\medskip
\begin{thmc}[van der Corput's Theorem]~\label{dth}
Let ${\bf x}=(x_{n})_{n\in \N}$ be a sequence of real numbers. If for every positive integer $h$ the sequence 
$$
{\bf d_{h}x} =(x_{n+h}-x_{n})_{n\in \N}
$$ 
is u.\,d.\,mod\;$1$, then ${\bf x}$ itself is u.\,d.\,mod\;$1$.
\end{thmc} 

The reader who is interested in these three theorems can find proofs in~\cite[Theorem 2.1, Theorem 3.1, and Theorem 4.3]{Kuipers-Niederreiter:74}.    

For the convenience of notation, we understand that an empty sum is $0$ 
and an empty product is one, i.e., $\sum_{j=1}^{k}(\cdots)=0$
and $\prod_{j=1}^k (\cdots)=1$ when $k=0$. In this notation, we have $y_{n} (x) =x^{n} =x^{n}\prod_{j=1}^{k} (\cdots)$ 
and $y_{n} (x) \sum_{j=1}^{} (\cdots) =0$ when $k=0$. 

\medskip
Take a non empty interval $I$ in the real line $\R$, which can be closed, open or semi-open.
Let $\mathcal{C}^{k}_{+}(I)$ be the space of all positive real valued $k$-times
continuously differentiable 
functions on an interval $I$, whose $i$-th derivative is non-negative for $i\le k$.
Then it is 
closed under addition and multiplication, that is, if $f,g\in \mathcal{C}^{k}_{+}(I)$, then $f+g, fg\in \mathcal{C}^{k}_{+}(I)$. In what follows, we often use this {\it closure property} of $\mathcal{C}^{k}_{+}(I)$. Let ${\mathbb R}[x]$ denote the space of all real polynomials.
The main result we prove in this paper is 

\medskip
\begin{thm}[Main Theorem]\label{mth1}
Let us take a function $g\in \mathcal{C}^{2}_{+}((1,\infty))$. Then, for a fixed real
number $\alpha\not=0$ and almost all real numbers $\beta>1$ 
(alternatively, for a fixed real number $\beta>1$ and almost all real numbers $\alpha$) and for all real polynomials $Q\in \R[x]$, 
sequences
\begin{equation}~\label{main1eq}
\Big( \alpha \beta^n g(\beta)+Q(n) \Big)_{n\in \N}
\end{equation}
are u.\,d.\,mod\;$1$.  
\end{thm}

As a countable union of exceptional null sets is null, we clearly have

\medskip
\begin{cor}\label{CSum}
Given a countable family $\{g_i\ | \ i\in \N \}$ in $\mathcal{C}^{2}_{+}((1,\infty))$.  Then, for a fixed real number $\alpha\not=0$ and almost all real numbers $\beta>1$ (alternatively, for a fixed real number $\beta>1$ and almost all real numbers $\alpha$)  and for all real polynomials $Q\in \R[x]$, sequences
$$
\Big( \alpha \beta^n g_i(\beta)+Q(n)\Big)_{n\in \N}, \quad \; i\in \N,
$$
are u.\,d.\,mod\;$1$.
\end{cor}

Since $x^n-1\in \mathcal{C}^{2}_{+}((1,\infty))$ for all $n\geq 1$,  
we have that for any $g\in \mathcal{C}^{2}_{+}((1,\infty))$ and any integer $l\geq 0$ and $h_{j} \in \N$, $1\leq j\leq l$,  
$$
g(x) \prod_{j=1}^{l} (x^{h_{j}}-1) \in \mathcal{C}^{2}_{+}((1,\infty)).
$$
Applying this countable family to Corollary \ref{CSum}, we obtain:

\medskip
\begin{thm}[Equivalent Statement of Main Theorem]~\label{mth2}
Given any $g\in \mathcal{C}^{2}_{+}((1,\infty))$. Then, for a fixed real number $\alpha\neq 0$ and almost all real numbers $\beta>1$
(alternatively, for a fixed real number $\beta>1$ and almost all real numbers $\alpha$), 
for all integers $l\geq 0$ and all $l$-tuple $(h_{1}, \cdots, h_{l})\in\N^{l}$, and for all real polynomials $Q\in \R[x]$, sequences
\begin{equation}~\label{main2eq}
\Big( \alpha \beta^n g(\beta) \prod_{j=1}^{l} (\beta^{h_{j}}-1)+Q(n)\Big)_{n\in \N}
\end{equation}
are u.\,d.\,mod\;$1$. 
\end{thm} 

On the other hand, when $l=0$, the product is empty, Theorem~\ref{mth2} is reduced to
Theorem~\ref{mth1}.  
So our Theorem~\ref{mth1} and Theorem~\ref{mth2} are equivalent.
Our proof is based on the formulation of Theorem~\ref{mth2}, which is already an interesting point in this paper. We will give the proof in \S3. 

Theorem~\ref{mth1} combined with Theorem A answers the question (Question~\ref{que}) affirmatively. 
 
\medskip
\begin{cor}[Main Corollary]~\label{mcoro}
Given any $g\in \mathcal{C}^{2}_{+}((1,\infty))$. Then, for a fixed real number $\alpha\not=0$ and almost all real numbers $\beta >1$ (alternatively, for a fixed real numbers $\beta>1$ 
and almost all real numbers $\alpha$), sequences
\begin{equation}~\label{mcor}
{\bf c} =\big( e^{2\pi i \alpha \beta^n  g(\beta)} \big)
\end{equation}
are higher order oscillating sequences of order $m$ for all $m\geq 2$ in the sense of the stronger 
version of the definition (Definition~\ref{hooss}).
\end{cor} 

\begin{rem}
In particular, taking a constant function $g\equiv 1$, we have ${\bf c}=(e^{2\pi i \alpha \beta^n})_{n\in \N}$.
\end{rem}

\section{Proof of the Main Theorem}

Our proof is based on the formulation of Theorem~\ref{mth2}. The reader may notices that it is already an interesting point in this paper. 
We start the proof for the 
case that $\alpha\neq 0$ is fixed and figure out the exceptional set for $\beta>1$.
We may assume that $\alpha>0$ and fixed it from now on. 
Take $g\in \mathcal{C}^{2}_{+}((1,\infty))$. 
For $n\in \N$ and an integer 
$l\geq 0$ and for a $l$-tuple $(h_1,\dots,h_l)\in \N^{l}$, define a function 
$$
y_n(x)=\alpha g(x) x^n \prod_{j=1}^l (x^{h_j}-1).
$$ 
Remember that when $l=0$, $y_{n}(x) =\alpha x^{n} g(x)$. 

For $n>m$, we have
\begin{eqnarray*}
y'_n(x)-y'_m(x)&=&
\alpha g(x)\Big((n x^{n-1}-m x^{m-1})\prod_{j=1}^k (x^{h_j}-1)\\
&+ &(x^n-x^m) \sum_{j=1}^k h_j x^{h_j-1} \prod_{i\neq j} (x^{h_i}-1)\Big)\\
&+&\alpha g'(x)(x^n-x^m)\prod_{j=1}^l (x^{h_j}-1).
\end{eqnarray*}
Since $g', g''\geq 0$ and since
$$
n x^{n-1}-m x^{m-1}=x^{m-1}(n x^{n-m}-m)
$$ 
and 
$$
n x^{n-m}-m\ge n-m\ge 1\quad \hbox{for $x\ge 1$},
$$ 
we see that every term in the last expression
are in $\mathcal{C}^{1}_{+}([1,\eta])$ for any $\eta>1$. 
By the closure property of $\mathcal{C}^{1}_{+}([1,\eta])$, we have that 
\begin{equation}~\label{inc}
y'_n(x)-y'_m(x)\in \mathcal{C}^{1}_{+}([1, \eta]) \quad \forall n>m\in \N. 
\end{equation}
In particular, this imply that $y'_n-y'_m$ is increasing for $n>m \in \N$.
Furthermore, if $x>a>1$, then $x^{h}-1\ge a-1$ for any $h\in \N$, we see that there is a constant $L>0$ such that
\begin{equation}~\label{min}
|y_n'(x)-y'_m(x)|\ge L\quad \forall n>m\in \N, \;\forall a \leq x\leq \eta. 
\end{equation}
Inequalities (\ref{inc}) and (\ref{min}) say that the sequence of real valued $C^{1}$ functions $(y_{n}(x))_{n\in \N}$ satisfies all hypothesizes of Theorem B. 

Theorem B implies that for almost all $x$ in $[(2^{k}+1)/2^{k}, (2^{k-1}+1)/2^{k-1}]$ or $[k, k+1]$ for $k\geq 2$, the sequence 
$(\alpha y_n(x))_{n\in \N}$ is u.\,d.\,mod\;$1$. Further, this implies that for almost all 
$$
x\in (1, \infty)= \cup_{k=2}^{\infty} \Big[\frac{2^{k}+1}{2^{k}}, \frac{2^{k-1}+1}{2^{k-1}}\Big]\cup  \cup_{k=2}^{\infty} [k, k+1]
$$ 
the sequence 
$(\alpha y_n(x))_{n\in \N}$ is u.\,d.\,mod\;$1$.

Let 
$$
A_{(h_1,h_2,\dots,h_l)}=\{ \beta >1 \ | \ (\alpha \beta^n g(\beta)
\prod_{j=1}^l (\beta^{h_j}-1))_{n\in \N} 
\text{ is not u.\,d.\,mod\;$1$ }  \}.
$$
Then it has one dimensional Lebesgue measure zero. By the above convention, we include the case $l=0$ as well. 

Since the set 
$$
U=\bigcup_{l=0}^{\infty} \{ (h_1,\dots,h_l)\ |\ h_j\in \N \}
$$ 
is countable, the one dimensional Lebesgue measure of
$$
\bigcup_{(h_1,\dots,h_l)\in U} A_{(h_1,\dots,h_l)} 
$$
is zero too.

For the fix a real number $\alpha\neq 0$ in the theorem,  take a real number 
$$
\beta\in (1,\infty) \setminus \bigcup_{(h_1,\dots,h_l)\in U} A_{(h_1,\dots,h_l)}.
$$
This says that the sequence $\big(\alpha \beta^n g(\beta)
\prod_{j=1}^l (\beta^{h_j}-1)\big)_{n\in \N} $ is u.\,d.\,mod\;$1$ for all integers $l\geq 0$ and all $l$-tuple $(h_{1}, \cdots, h_{l})\in \N^{l}$.

Define statements $P(k)$ for $k=0,1,\dots$ as follows. 

\medskip
{\sl
$P(k)$: For any integer $l\geq 0$, 
$(h_1,\dots,h_{l})\in \N^{l}$ and $t_i\in \R\ (i=0,1,2,\dots, k)$, the sequence 
$$
\left(\alpha \beta^n g(\beta) \prod_{j=1}^{l} (\beta^{h_j}-1)+ \sum_{i=0}^k t_i n^i \right)_{n\in \N}
$$
is u.\,d.\,mod\;$1$.
}

\medskip
We claim that $P(k)$ holds for every integer $k\geq 0$. We prove the claim by induction. 

By our choice of $\alpha$ and $\beta$, we know that $P(0)$ holds. 

Assume $P(k-1)$ holds for $k\geq 1$. Let
$$
x_{n}=\alpha \beta^n g(\beta)\prod_{j=1}^{\ell} (\beta^{h_j}-1)+ \sum_{i=0}^k t_i n^i.
$$
Then
$$
x_{n+h}=\alpha \beta^{n+h} g(\beta)
\prod_{j=1}^{\ell} (\beta^{h_j}-1)+ \sum_{i=0}^k t_i (n+h)^i.
$$
Consider the difference appeared in Theorem C:
\begin{eqnarray*}
x_{n+h}-x_{n}&=&\alpha \beta^n (\beta^h-1) g(\beta)\prod_{j=1}^{\ell} (\beta^{h_j}-1)+ \sum_{i=0}^{k-1}T_i n^i\\
&=&\alpha \beta^n g(\beta)\prod_{j=1}^{\ell +1} (\beta^{h_j}-1)+ \sum_{i=0}^{k-1}T_i n^i\end{eqnarray*}
with
$$
h_{\ell +1}=h \hbox{ and } T_i=-t_i+\sum_{j=i}^k t_j {j \choose i} h^{j-i}.
$$
Since $P(k-1)$ is valid, the resulting sequence is u.\,d.\,mod\;$1$ for all $h\in \N$. Now Theorem C implies that $P(k)$
holds too. We proved the claim.  Therefore we completed the proof of Theorem~\ref{mth2}. 
  
The proof for the 
case that $\beta>1$ is fixed and to obtain the exceptional set for $\alpha$, is similar and easier. 
Let $g$ be any positive function on $(1,\infty)$.
For $n\in \N$ and an integer $l\geq 0$ and for a $l$-tuple $(h_1,\dots,h_l)\in \N^{l}$, define 
$$
y_n(x)=x g(\beta) \beta^n \prod_{j=1}^l (\beta^{h_j}-1).
$$ 
Then $y'_n(x)-y'_m(x)=g(\beta)(\beta^n-\beta^m)\prod_{j=1}^l (\beta^{h_j}-1)$ 
is a positive constant for $n>m$ and 
satisfies
the condition of Theorem B under the similar dissection of the interval $(1,\infty)$.
The rest of the proof is the same.

\medskip
\medskip
\noindent {\em Acknowledgment.} This work was done when both of the authors visited the National Center for Theoretical Sciences (NCTS) at National Taiwan University during 2016. They would like to thank NCTS for its hospitality. We also like to thank Professors Jung-Chao Ban and Chih-Hung Chang and other audiences for their spending times patiently to listen and discuss lectures given by both of the authors in NCTS including proofs in this paper.

\medskip
\medskip

\noindent Shigeki AKIYAMA:
Institute of Mathematics, University of Tsukuba
1-1-1 Tennodai, Tsukuba, Ibaraki, 305-8571 Japan.\\
Email:akiyama@math.tsukuba.ac.jp 

\medskip
\noindent Yunping JIANG: Department of Mathematics, 
Queens College of the City University of New York,
Flushing, NY 11367-1597 USA and 
Department of Mathematics
Graduate School of the City University of New York,
365 Fifth Avenue, New York, NY 10016 USA.\\
Email:yunping.jiang@qc.cuny.edu

\end{document}